\documentclass{article}
\usepackage{latexsym}
\usepackage{amsmath}
\usepackage{eufrak}

\catcode`\@=11
\def\section{\@startsection{section}{2}{\z@}{-3.25ex plus -1ex minus 
 -.2ex}{1.5ex plus .2ex}{\large\bf}}
\catcode`\@=12
\newcounter{hilf2}
\def\listzero{\topsep0pt
              \parskip0pt
              \parsep0pt
              \itemsep 1.5ex plus 0.5ex minus 0.5ex
              \leftmargin0pt
              \listparindent0pt
              \labelwidth20pt
              \labelsep5pt
              \itemindent25pt}

\def\myenumerate#1{\list{(\roman{hilf2})}{\usecounter{hilf2}\listzero #1}}
\newenvironment{enumerateremark}{%
\myenumerate{\listparindent\parindent}}{\endlist}
\catcode`\@=11
\newcommand{\model}{\mathit{M}}
\newcommand{\theory}{\mathit{T}}
\newtheorem{theorem}{Theorem}[section]
\newtheorem{lemma}[theorem]{Lemma}

\newtheorem{proposition}[theorem]{Proposition}
\newtheorem{definition}[theorem]{Definition}
\newtheorem{remark}[theorem]{Remark}
\newenvironment{proof}{\trivlist\item[]\textit{Proof.\/}
                       }{\hspace*{\fill}$\Box$\endtrivlist}
\newenvironment{proof*}[1]{\trivlist\item[]\textit{#1\/}}{\hspace*{\fill}$\Box$\endtrivlist}
\newcommand{\logic}{\mathcal{L}}
\newcommand{\modM}{\mathfrak{M}}
\newcommand{\card}[1]{|#1|}
\newcommand{\valuate}[1]{[\![ #1 ]\!]}
\newcommand{\equal}[2]{\valuate{#1 = #2}}
\newcommand{\BA}{B}
\newcommand{\points}{\mathcal{S}_{\theory}}
\newcommand{\nodel}{\mathrm{N}}

\catcode`\@=12
\begin{document}
\title{A definability theorem for first order logic
       }
\author{Carsten Butz ({\AA}rhus) and Ieke Moerdijk (Utrecht)\thanks{%
Both authors acknowledge support from the 
Netherlands Science Organisation~(NWO).}
        }
\date{\relax}
\maketitle
%
%
%
%

In this paper, we will present a definability theorem for first order logic. 
This theorem is very easy to state, and its proof only uses elementary tools. 
To explain the theorem, let us first observe that if $\model$ is a model of 
a theory $\theory$ in a language~$\logic$, then, clearly, any definable subset 
$S\subset \model$ (i.e.,~a subset $S=\{a\mid \model\models\varphi(a)\}$ 
defined by some formula~$\varphi$) is invariant under all automorphisms 
of~$\model$. The same is of course true for subsets of $\model^n$ defined by 
formulas with $n$ free variables.

Our theorem states that, if one allows Boolean valued models, the converse holds. More precisely, for any theory $\theory$ we will construct a Boolean valued model~$\model$, in which precisely the $\theory$--provable formulas hold, and in which every (Boolean valued) subset which is invariant under all automorphisms of $\model$ is definable by a formula of~$\logic$.

Our presentation is entirely selfcontained, and only requires familiarity with the most elementary properties of model theory. In particular, we have added a first section in which we review the basic definitions concerning Boolean valued models.

The Boolean algebra used in the construction of the model will be presented concretely as the algebra of closed and open subsets of a topological space $X$ naturally associated with the theory~$\theory$.
The construction of this space is closely related to the one in~\cite{Butz-Moerdijk:96c}.
In fact, one of the results in that paper could be interpreted as a definability theorem for infinitary logic, using topological rather than Boolean valued models.

\section{Preliminary definitions}
\label{preliminary}

In this section we review the basic definitions concerning Boolean valued models
(see e.g.~\cite{Koppelberg:85}).
Most readers will be familiar with these notions, and they are advised to skip this section. They should note, however, that our Boolean algebras are not necessarily complete, and that we treat constants and function symbols as functional relations.

Let us fix a signature~$S$, consisting of constants, function and relation symbols. For simplicity we assume it is a single sorted signature, although this restriction is by no means essential. Let $\logic$ denote the associated first order language $\logic_{\omega\omega}(S)$.

A {\em Boolean valued interpretation\/} of $\logic$ is a triple
$\modM=(\BA,\card{\modM},\valuate{-})$,
where $\BA$ is a Boolean algebra,
$\card{\modM}$ is the underlying set of the interpretation,
and $\valuate{-}$ is an operation which assigns to each formula 
$\varphi(x_1,\ldots,x_n)$ of $\logic$ with free variables among 
$x_1,\ldots,x_n$ a function
$\card{\modM}^n\to\BA$, whose value at $(m_1,\ldots,m_n)$ is denoted
$$\valuate{\varphi(m_1,\ldots,m_n)}.$$
These functions are required to satisfy the usual identities 
(where we write $m$ for $m_1,\ldots,m_n$):
\begin{enumerate}
\item[(i)] %
     $\valuate{\varphi\wedge\psi(m)}=
                 \valuate{\varphi(m)}\wedge\valuate{\psi(m)}$
     and similar for the other Boolean connectives.
\item[(ii)] %
     $\valuate{\exists y\varphi(y,m)}=
             \bigvee\{\valuate{\varphi(k,m)}\mid k\in\card{\modM}\}$,\\[1ex]
     $\valuate{\forall y\varphi(y,m)}=
             \bigwedge\{\valuate{\varphi(k,m)}\mid k\in\card{\modM}\}$,
\end{enumerate}
 where it is part of the definition of an interpretation that these sups and infs are required to exists in~$B$.
Finally, we require
\begin{enumerate}
\item[(iii)] %
     if $\vdash\varphi(x_1,\ldots,x_n)$
     then $\valuate{\varphi(m)}=1_{\BA}$
     for any $m\in\card{\modM}^n$.
\end{enumerate}
In (iii), $\vdash$ denotes derivability in (one of the usual axiomatisations of) classical first order logic.

\begin{remark}
\begin{rm}
\begin{enumerateremark}
\item \leftmargin0pt
Note that, in particular, $\card{\modM}$ is equipped with a $\BA$--valued equality $\equal{x_1}{x_2}\colon\card{\modM}^2\to\BA$,
satisfying the identities for reflexivity, transitivity and symmetry,
$$\begin{array}{l}
\equal{m}{m}=1_\BA,\\[1ex]
\equal{m_1}{m_2}=\equal{m_2}{m_1},\\[1ex]
\equal{m_1}{m_2}\wedge\equal{m_2}{m_3}\leq\equal{m_1}{m_3}.
\end{array}
$$
\item For each constant $c$ the formulas $c=x$ and $x=c$ define the same function $C=\equal{c}{x}\colon\card{\modM}\to\BA$,
which should be viewed as the interpretation of~$c$. It satisfies the conditions
$C(m)\wedge\equal{m}{m'}\leq C(m')$
and
$\bigvee\{C(m)\mid m\in\card{\modM}\}=1_\BA$.
Similarly, each $n$--ary function symbol is interpreted, via the formula 
$f(x_1,\ldots,x_n)=y$, 
by a function $F\colon\card{\modM}^n\times\card{\modM}\to\BA$.
This function satisfies the conditions
$F(m,k)\wedge\equal{m}{m'}\wedge\equal{k}{k'}\leq F(m',k')$
and
$\bigvee\{F(m,k)\mid k\in \card{\modM}\}=1_\BA$.
(Here $m=m_1,\ldots,m_n$ as before, and $\equal{m}{m'}$ stands for
$\bigwedge_{i=1}^n\equal{m_i}{m_i'}$.)
\item For each $n$--ary relation symbol $r$ the formula $r(x_1,\ldots,x_n)$ 
defines a map
$R\colon\card{\modM}^n\to\BA$, which is extensional in the sense that
$R(m)\wedge\equal{m}{m'}\leq R(m')$.
\item The entire interpretation is determined by these data in (i)--(iii).
First, using derivability of usual equivalences such as
$\vdash f(g(x))=y\leftrightarrow\exists z(f(z)=y\wedge g(x)=z)$,
one obtains by induction for each term $t(x_1,\ldots, x_n)$ a function 
$T\colon\card{\modM}^{n+1}\to\BA$ interpreting the formula 
$t(x_1,\ldots,x_n)=y$.
Next, one builds up the interpretation of formulas in the usual way, using the assumption that all necessary sups and infs exist in~$\BA$.
\end{enumerateremark}
\end{rm}
\end{remark}

As usual, we write $\modM\models\varphi$ if $\valuate{\varphi(m)}=1$ for all 
$m\in\card{\modM}^n$, and we say $\modM$ is a model of a theory $\theory$ if $\modM\models\varphi$ whenever~$\theory\vdash\varphi$. In this case, we write $\modM\models\theory$, as usual.

\section{Automorphisms of models and statement of the theorem}

Consider a fixed Boolean valued model $\modM=(\BA,\card{\modM},\valuate{-})$.
An {\em automorphism\/} $\pi$ of $\modM$ consists of two mappings $\pi_0$ and~$\pi_1$. The map $\pi_0\colon\BA\to\BA$ is an automorphism of the Boolean algebra~$\BA$, while $\pi_1\colon\card{\modM}\to\card{\modM}$ is an automorphism of the underlying set~$\card{\modM}$, with the property that 
\begin{equation}\label{eq:2.1}
\pi_0\valuate{\varphi(m_1,\ldots,m_n)}=
  \valuate{\varphi(\pi_1(m_1),\ldots,\pi_1(m_n))},
\end{equation}
for any formula $\varphi(x_1,\ldots,x_n)$ and any $m_1,\ldots,m_n\in\card{\modM}$.
(Of course it is enough to check a condition like~(\ref{eq:2.1}) for constants, functions and relations of~$\logic$, and deduce~(\ref{eq:2.1}) for arbitrary $\varphi$ by induction.)

An ($n$--ary) {\em predicate\/} on $\modM$ is a map $p\colon\card{\modM}^n\to\BA$ which satisfies the extensionality condition
\begin{equation}\label{eq:2.2}
p(m)\wedge\equal{m}{m'}\leq p(m')
\end{equation}
for any $m,m'\in\card{\modM}^n$
(where $\equal{m}{m'}$ stands for
$\bigwedge_{i=1}^n\equal{m_i}{m'_i}$, as before).
Such a predicate $p$ is {\em definable\/} if there is a formula $\varphi(x_1,\ldots,x_n)$ such that
\begin{equation}\label{eq:2.3}
p(m)=\valuate{\varphi(m)},\qquad\mbox{for all $m\in\card{\modM}^n$.}
\end{equation}
It is {\em invariant\/} under an automorphism $\pi$ if
\begin{equation}\label{eq:2.4}
\pi_0p(m)=p(\pi_1(m)),\qquad\mbox{for all $m\in\card{\modM}^n$,}
\end{equation}
(where $\pi_1(m)$ is $(\pi_1(m_1),\ldots,\pi_1(m_n))$).
Obviously, every definable predicate is invariant. Our theorem states the converse.

\begin{theorem}\label{main-theorem}
Let $\theory$ be any first order theory. There exists a Boolean valued model 
$\modM$ such that
\begin{enumerate}
\item $\modM$ is a conservative model of~$\theory$, in the sense that 
      $\modM\models\varphi$ iff $\theory\vdash\varphi$, 
      for any sentence~$\varphi$.
\item Any predicate which is invariant under all automorphisms of $\modM$ is 
      definable.
\end{enumerate}
\end{theorem}

Before proving the theorem in \S\ref{proof}, we will first give an explicit 
description of the Boolean algebra and the interpretation involved in the next 
section.

\section{Construction of the model}
\label{the-model}

Our Boolean algebra will be defined as the algebra of all clopen 
(i.e.,~closed and open) sets in a topological space~$X$.
To describe~$X$, let $\kappa\geq\omega$ be the cardinality of our 
language~$\logic$.
We fix a {\em set\/} $\points$ of (ordinary, two--valued) models $\model$ of 
$\theory$ such that every model of cardinality $\leq\kappa$ is isomorphic to 
a model in~$\points$.
Then, in particular, a formula is provable from $\theory$ iff it holds in all 
models in the set~$\points$.

\begin{definition}
An {\em enumeration\/} of a model $\model$ is a function $\alpha\colon\kappa\to\card{\model}$ such that $\alpha^{-1}(a)$ is infinite for all $a\in\card{\model}$
(here $\card{\model}$ is the underlying set of~$\model$).
\end{definition}

The space $X$ has as its points the equivalence classes of pairs $(\model,\alpha)$, where $\model\in\points$ and $\alpha$ is an enumeration of~$\model$.
Two such pairs $(\model,\alpha)$ and $(\nodel,\beta)$ are {\em equivalent\/} if there exists an isomorphism of models $\theta\colon\model\stackrel{\simeq}{\to}\nodel$ such that $\beta=\theta\circ\alpha$.
We will often simply write $(\model,\alpha)$ when we mean the equivalence class of~$(\model,\alpha)$.
The topology of $X$ is generated by all the basic open sets of the form
\begin{equation}\label{3.1}
U_{\varphi,\xi}=\{(\model,\alpha)\mid\model\models\varphi(\alpha(\xi))\}.
\end{equation}
Here $\varphi=\varphi(x_1,\ldots,x_n)$ is any formula with free variables among $x_1,\ldots,x_n$, while $\xi=(\xi_1,\ldots,\xi_n)$ is a sequence of elements of~$\kappa$ (i.e.,~ordinals $\xi_i<\kappa$);
we use $\alpha(\xi)$ as an abbreviation of $\alpha(\xi_1),\ldots,\alpha(\xi_n)$.

Observe that each such basic open set $U_{\varphi,\xi}$ is also closed, with complement $U_{\neg\varphi,\xi}$. So $X$ is a zero--dimensional space. 
We now define the Boolean algebra $\BA$ as
\begin{equation}\label{eq:3.2}
\BA=\mathrm{Clopens}(X),
\end{equation}
the algebra of all open and closed sets in~$X$.

Notice that arbitrary suprema need not exist in~$\BA$, although $\BA$ has many infinite suprema.
In particular, if $U\subset X$ is clopen and $\{U_i\}_{i\in I}$ is a cover of $U$ by basic open sets, then the union $\bigcup_{i\in I}U_i$ defines a supremum $U=\bigvee_{i\in I}U_i$ in~$\BA$; we only need suprema of this kind.
\vskip 2ex

The Boolean algebra $\BA$ just constructed is part of a natural Boolean valued model $\modM=(\BA,\card{\modM},\valuate{-})$, with
\begin{equation}\label{eq:4.1}
\card{\modM}=\kappa
\end{equation}
and evaluation of formulas defined by
\begin{equation}\label{eq:4.2}
\valuate{\varphi(\xi_1,\ldots,\xi_n)}=U_{\varphi,\xi},
\end{equation}
for any formula $\varphi(x_1,\ldots,x_n)$ and any sequence 
$\xi=\xi_1,\ldots,\xi_n$ of ordinals~$\xi_i<\kappa$.

\begin{lemma}\label{4.1}
This evaluation defines a $\BA$--valued interpretation of the language~$\logic$.
\end{lemma}

\begin{proof}
One needs to check the requirements (i)--(iii) from Section~\ref{preliminary}.
Now~(iii) is clear, while (i) and (ii) are completely straightforward.
For illustration, we give the case of the existential quantifier.
Suppose $\varphi(y,x)$ is a formula with just two free variables $x$ and $y$.
Then for any $\xi<\kappa$,
\begin{eqnarray*}
\valuate{\exists y\varphi(y,\xi)}
& = & \{(\model,\alpha)\mid \model\models\exists y\varphi(y,\alpha(\xi))\}\\
& = & \{(\model,\alpha)\mid \exists\eta<\kappa\colon\
             \model\models\varphi(\alpha(\eta),\alpha(\xi))\}
\\
& & \mbox{(since each $\alpha$ is surjective)}\\
& = & \bigcup_{\eta<\kappa}\{(\model,\alpha)\mid
             \model\models\varphi(\alpha(\eta),\alpha(\xi))\}\\
& = & \bigcup_{\eta<\kappa}\valuate{\varphi(\eta,\xi)},
\end{eqnarray*}
and this union is a supremum in~$\BA$, by the remark above.
\end{proof}

\section{Proof of the theorem}
\label{proof}

We will now show that the interpretation $\modM$ has the two properties stated 
in Theorem~\ref{main-theorem}. The first one is easy:

\begin{proposition}
$\modM$ is a conservative model of~$\theory$.
\end{proposition}

\begin{proof}
We need to show that $\modM\models\sigma$ iff $\theory\vdash\sigma$, for any sentence $\sigma\in\logic$.
By Lemma~\ref{4.1}, 
$\valuate{\sigma}=\{(\model,\alpha)\mid\model\models\sigma\}$.
Thus $\valuate{\sigma}=X$ iff $\model\models\sigma$ for all $\model\in\points$, and this holds iff $\theory\vdash\sigma$, by definition of~$\points$.
\end{proof}

For the proof of the definability result~\ref{main-theorem}(ii),
we shall only need a particular collection of automorphisms of the model~$\modM$.
Let $S_\kappa$ denote the symmetric group of permutations of~$\kappa$. Then $S_\kappa$ acts on the model $\modM$ as follows. Any $\pi_1\in S_\kappa$ induces a homeomorphism $\pi_0\colon X\to X$, defined by
$$\pi_0(\model,\alpha)=(\model,\alpha\circ\pi_1^{-1}).$$
This map has the property that 
$\pi_0(U_{\varphi,\xi})=U_{\varphi,\pi_1(\xi)}$, or
$$\pi_0\valuate{\varphi(\xi)}=\valuate{\varphi(\pi_1(\xi))},
$$
for any formula $\varphi(x_1,\ldots,x_n)$ and any $\xi=\xi_1,\ldots,\xi_n<\kappa$.
Thus, the pair $\pi=(\pi_1,\pi_0)$ is an automorphism of~$\modM$. This defines an action of $S_\kappa$ on~$\modM$, i.e., a representation
$$\rho\colon S_\kappa\to\mathrm{Aut}(\modM),\qquad\rho(\pi_1)=\pi.
$$
For the second part of Theorem~\ref{main-theorem}, it will now be enough to show:

\begin{proposition}\label{5.2}
Any $S_\kappa$--invariant predicate is definable.
\end{proposition}

To simplify notation, we will only prove this for a unary predicate.
So let us fix such an invariant predicate~$p$.
It is a function
$p\colon\card{\modM}=\kappa\longrightarrow\BA$ satisfying the extensionality 
condition
$$p(\xi)\wedge\equal{\xi}{\xi'}\leq p(\xi'),$$
as well as the invariance condition
$$ p(\pi_1\xi)=\pi_0(p(\xi)),$$
for any $\pi_1\in S_\kappa$.
We will first show that $p$ is ``locally'' definable (Lemma~\ref{5.5}).

\begin{lemma}\label{5.3}
Let $(\model,\alpha)\in U\in \BA$ and $\eta_0\in\kappa$.
Then there is a formula
$\delta(x_1,\ldots,x_n,y)$ and elements $\xi_1,\ldots,\xi_n\in\kappa$ such that
\begin{enumerate}
\item $(\model,\alpha)\in U_{\delta,(\xi,\eta_0)}\leq U$.
\item For any point $(\nodel,\beta)$ in $X$, any $b_1,\ldots,b_n,c\in\card{\nodel}$ such that
$\nodel\models\delta(b_1,\ldots,b_n,c)$, and any $\eta\in\kappa$ with $\beta(\eta)=c$,
there exists a $\pi_1\in S_\kappa$ such that $\pi_1(\eta)=\eta_0$ and
$\pi_0(\nodel,\beta)\in U_{\delta,(\xi,\eta_0)}$.
\end{enumerate}
\end{lemma}

\begin{proof}
Choose a basic open set $U_{\delta',\xi}$, given by a formula 
$\delta'(x_1,\ldots,x_n)$ and $\xi_1,\ldots,\xi_n<\kappa$, such that
$$(\model,\alpha)\in U_{\delta',\xi}\subset U.$$
Let $\mathrm{Eq}_\alpha(x_1,\ldots,x_n,y)$ be the formula
$$\bigwedge_{\alpha(\xi_i)=\alpha(\xi_j)}x_i=x_j\wedge
   \bigwedge_{\alpha(\xi_i)=\alpha(\eta_0)}x_i=y,
$$
and define $\delta$ to be $\delta'\wedge\mathrm{Eq}_\alpha$. Then obviously
$$(\model,\alpha)\in U_{\delta,\xi,\eta_0}\subset U_{\delta',\xi}\subset U.
$$
Now choose any $(\nodel,\beta)$, $b_1,\ldots,b_n,c$ and $\eta$ satisfying the hypothesis of part~(ii) of the lemma. Then in particular
$\nodel\models\mathrm{Eq}_\alpha(b_1,\ldots,b_n,c)$ and $c=\beta(\eta)$. Since $\beta\colon\kappa\to\card{\nodel}$ has infinite fibres, we can find $\zeta_1,\ldots,\zeta_n<\kappa$ such that $\beta(\zeta_i)=b_i$, while the sequence $\zeta_1,\ldots,\zeta_n,\eta$ satisfies {\em exactly\/} the same equalities and inequalities as the sequence $\xi_1,\ldots,\xi_n,\eta_0$.
[Indeed, if $\zeta_1,\ldots,\zeta_i$ have been found, and $\xi_{i+1}=\xi_k$ for some $k\leq i$ or $\xi_{i+1}=\eta_0$, then also 
$\alpha(\xi_{i+1})=\alpha(\xi_k)$ or $\alpha(\xi_{i+1})=\alpha(\eta_0)$, hence $b_{i+1}=b_k$ or $b_{i+1}=c$ since $\nodel\models \mathrm{Eq}_\alpha(b_1,\ldots,b_n,c)$. Thus, we can choose $\zeta_{i+1}=\zeta_k$ respectively $\zeta_{i+1}=\eta$.
If, on the other hand, $\xi_{i+1}\notin\{\eta_0,\xi_1,\ldots,\xi_i\}$, we can use the fact that $\beta^{-1}(b_{i+1})$ is infinite, to find $\zeta_{i+1}\in\beta^{-1}(b_{i+1})\setminus\{\eta,\zeta_1,\ldots,\zeta_i\}$.]
Thus, there is a permutation $\pi_1\in S_\kappa$ with 
$$\pi_1(\eta)=\eta_0,
\pi_1(\zeta_1)=\xi_1,\ldots,
\pi_1(\zeta_n)=\xi_n.$$
But then $\nodel\models\delta(b_1,\ldots,b_n,c)$ means that 
$\nodel\models\delta(\beta(\pi_1^{-1}(\xi_1)),\ldots,
                     \beta(\pi_1^{-1}(\xi_n)),$ $\beta(\pi_1^{-1}(\eta_0)))$,
or that $\pi_0(\nodel,\beta)\in U_{\delta,(\xi,\eta_0)}$.
\end{proof}

\begin{lemma}\label{5.4}
Let $\eta_0<\kappa$. There is a cover $p(\eta_0)=\bigvee_{i \in I(\eta_0)}U_i$ in~$\BA$, and formulas $\psi_i^{\eta_0}(y)$, such that for any $i\in I(\eta_0)$,
\begin{enumerate}
\item $U_i\leq\valuate{\psi_i^{\eta_0}(\eta_0)}$.
\item For any $\eta<\kappa$, $\valuate{\psi_i^{\eta_0}(\eta)}\leq p(\eta)$.
\item$\bigvee\limits_{i\in I(\eta_0)}\valuate{\psi_i^{\eta_0}(\eta_0)}=p(\eta_0)$.
\end{enumerate}
\end{lemma}

\begin{proof} Observe that (iii) follows from (i) and~(ii). To prove these, write $U=p(\eta_0)$, and apply Lemma~\ref{5.3} to each of the points $(\model,\alpha)\in U$. This will give a cover $U=\bigcup_{i\in I}U_i$ by basic open sets, and for each index $i$ a formula $\delta_i(x_1,\ldots,x_n,y)$ and elements $\xi_1,\ldots,\xi_n<\kappa$ such that
$$U_i=U_{\delta_i,(\xi,\eta_0)},$$
and moreover such that property~(ii) of Lemma~\ref{5.3} holds for each of these formulas~$\delta_i$. Now define
$$\psi_i^{\eta_0}(y)=\exists x_1\ldots\exists x_n\delta_i(x_1,\ldots,x_n,y).$$
It is now clear that statement~(i) in the lemma holds. For~(ii), suppose $(\nodel,\beta)\in\valuate{\psi_i^{\eta_0}(\eta)}$.
This means that
$\nodel\models \exists x_1\ldots\exists x_n\delta_i(x_1,\ldots,x_n,\beta(\eta))$.
By~\ref{5.3}(ii), we can find a $\pi_1\in S_\kappa$ such that $\pi_1(\eta)=\eta_0$ and $\pi_0(\nodel,\beta)\in U_{\delta_i,(\xi,\eta_0)}=U_i$.
Since
$U_i\subset U=p(\eta_0)$, also $\pi_0(\nodel,\beta)\in p(\eta_0)$, and hence, by invariance of~$p$,
$(\nodel,\beta)\in p(\pi_1^{-1}(\eta_0))=p(\eta)$, as required.
\end{proof}

\begin{lemma}\label{5.5}
There is a family $\{\psi_i(y)\mid i\in I\}$ of formulas such that, for all $\eta<\kappa$,
$$p(\eta)=\bigvee_{i\in I}\valuate{\psi_i(\eta)}.$$
\end{lemma}

\begin{proof}
This follows immediately from the previous lemma, for the collection of formulas
$\{\psi_i^{\eta_0}\mid\eta_0<\kappa,\ i\in I(\eta_0)\}$.
\end{proof}

\begin{proof*}{Proof of Proposition~\ref{5.2}.}
Consider the function $p'\colon\card{\modM}\to\BA$ defined by $p'(\eta)=\neg p(\eta)$. Clearly, since $p$ is a predicate, so is~$p'$, i.e., $p'(\eta)\wedge\equal{\eta}{\eta'}\leq p'(\eta')$ for all $\eta,\eta'<\kappa$.
Moreover, $p'$~is invariant since $p$ is. So we can apply Lemma~\ref{5.5} to~$p'$, to find a collection of formulas
$$
\{\varphi_j(y)\mid j\in J\}
$$
such that for all $\eta<\kappa$,
\begin{equation}\label{eq:5.1}
p'(\eta)=\bigvee_{j\in J}\valuate{\varphi_j(\eta)}.
\end{equation}

The definability of $p$ now follows by a standard compactness argument. 
Let $c$ be a ``new'' constant, and consider the theory 
$\theory'=\theory\cup\{\neg\psi_i(c)\mid i\in I\}
                  \cup\{\neg\varphi_j(c)\mid j\in J\}$. 
If $\theory'$ where consistent, it would have a model~$\model$, which we can 
assume to be (an expansion of a model) in the set~$\points$. 
Let $\alpha$ be an enumeration of~$\model$, and 
choose $\eta<\kappa$ with $\alpha(\eta)=c^{(\model)}$, the interpretation of 
$c$ in~$\model$. 
Then $(\model,\alpha)\in X=p(\eta)\vee p'(\eta)$, hence 
$(\model,\alpha)\in\valuate{\psi_i(\eta)}$ for some $i\in I$ or 
$(\model,\alpha)\in\valuate{\varphi_j(\eta)}$ for some $j\in J$.
This means that 
$\model\models\psi_i(\alpha(\eta))\vee\varphi_j(\alpha(\eta))$, 
contradicting the fact that $\model$ models~$\theory'$.
This proves that $\theory'$ is inconsistent. 

Now apply compactness, to find $i_1,\ldots,i_n\in I$ and $j_1,\ldots,j_m\in J$ 
such that
\begin{equation}\label{eq:5.2}
\theory\vdash\forall y(\psi_{i_1}(y)\vee\cdots\vee\psi_{i_n}(y)\vee
                       \varphi_{j_1}(y)\vee\cdots\vee\varphi_{j_m}(y)).
\end{equation}
Write $\psi=\psi_{i_1}\vee\cdots\vee\psi_{i_n}$ and
$\varphi=\varphi_{j_1}\vee\cdots\vee\varphi_{j_m}$.
We claim that $\psi$ defines~$p$. Indeed, let $(\model,\alpha)$ be any point 
in~$X$, and let $\eta<\kappa$.
By~(\ref{eq:5.2}), $\model\models\psi(\alpha(\eta))\vee\varphi(\alpha(\eta))$,
or in other words, either $(\model,\alpha)\in\valuate{\psi(\eta)}$ or 
$(\model,\alpha)\in\valuate{\varphi(\eta)}$.
If $(\model,\alpha)\in\valuate{\psi(\eta)}$, then $(\model,\alpha)\in p(\eta)$ 
by Lemma~\ref{5.2}.
And if $(\model,\alpha)\in \valuate{\varphi(\eta)}$, then 
$(\model,\alpha)\in p'(\eta)$ by~(\ref{eq:5.1}), 
hence $(\model,\alpha)\notin p(\eta)$. 
Thus $(\model,\alpha)\in \valuate{\psi(\eta)}$ iff 
$(\model,\alpha)\in p(\eta)$.

This shows that $\valuate{\psi(\eta)}=p(\eta)$ for any $\eta<\kappa$, 
and completes the proof.
\end{proof*}

%
%
%
%

%
%
\vskip 3ex
\begin{minipage}{\textwidth}
\footnotesize
Ieke Moerdijk,
Mathematisch Instituut,
Universiteit Utrecht,
Postbus 80.010,
NL--3508 TA Utrecht,
The Netherlands,
moerdijk@math.ruu.nl.\\[1ex]
Carsten Butz,
BRICS, 
Computer Science Department, 
Aarhus University, 
Ny Munkegade, Building 540, DK-8000 {\AA}rhus~C, Denmark,
butz@brics.dk.\\[1ex]
{\bf BRICS}\\
Basic Research in Computer Science,
Centre of the Danish National Research Foundation.
\end{minipage}

\begin{thebibliography}{1}


\bibitem{Butz-Moerdijk:96c}
C.~Butz and I.~Moerdijk.
\newblock Representing topoi by topological groupoids.
\newblock Utrecht preprint nr.~984, October 1996. 
          To appear in J.~Pure and Applied Alg.

\bibitem{Koppelberg:85}
S.~Koppelberg.
\newblock Booleschewertige Logik.
\newblock Jber.\ d.~dt. Math.-Verein~\textbf{87} (1985), 19--38.

\end{thebibliography}
\end{document}